\documentstyle[leqno,12pt,twoside]{article}
\begin{document}

\begin{titlepage}

\title{\large\bf On Higher Regularized Traces of a Differential Operator with Bounded
Operator Coefficient Given in a Finite Interval}

\pagenumbering{arabic}\setcounter{page}{1}\pagestyle{myheadings}
\baselineskip=18pt \maketitle

\centerline {\bf Serpil Karayel, Yonca Sezer and \"{O}zlem Bak\c{s}i }

\centerline{Department of Mathematics,}

\centerline{ Faculty of Arts and Science, Y\i ld\i z Technical University }
\centerline {(34210), Davutpa\d{s}a, \.Istanbul, Turkey}

\vskip 0.3cm

\noindent
\begin{abstract} In this work, a higher regularized trace formula has been found for
a regular Sturm-Liouville differential operator with operator coefficient.

\end{abstract}

\noindent {\small  {\bf AMS Subj. Classification:} 47A10, 34L20,
34L05

}
\end{titlepage}


\section{Introduction}
Let $\>H\>$ be an infinite-dimensional separable Hilbert space. We will denote by
$\>(.,.)_{H}\>$ inner product and denote by $\>\|.\|_{H}\>$ the norm in $\>H\>$. Let
$\>H_{1}=L_{2}(0,\pi;\>H)\>$ be the set of all strongly measurable functions $\>f\>$ defined
on $\>[0,\pi]\>$ and taking values in the space $\>H\>$ satisfying the following conditions:

\newcommand{\ud}{\mathrm {d}}

\noindent 1. The scalar function $\>\bigl(f(x),g \bigr)_{H}\>$ is Lebesgue
measurable on the interval $\>[0,\pi]\>$, for every $\>g\in H$.\\
\\
\noindent 2. $\int\limits_{0}^{\pi}\|f(x)\|_{H}^{2} \>\ud x<\infty$.

\noindent  $\>H_{1}\>$ is a linear space. We will denote by $\>(.,.)\>$ the inner product
 in $\>H_{1}\>$ and by $\>\|.\|\>$ the norm in $\>H_{1}\>$. If the inner product is defined by

\begin{displaymath}
(f_{1},f_{2})=\int\limits_{0}^{\pi}\Bigl(f_{1}(x),f_{2}(x)\Bigr)_{H}\>\ud x
\end{displaymath}

\noindent for two arbitrary elements $\>f_{1} \>\textrm{and}\> f_{2}\>$ of $\>H_{1}\>$
then the space $\>H_{1}\>$ becomes a separable Hilbert
space,
\cite{Ki}.

\noindent Let us consider the operators $\>L_{0}\>$ and $\>L\>$ formed by the following
differential expressions
\begin{eqnarray*}
\qquad\qquad\qquad\ell_{0}(y)&=&-y''(x)\quad \qquad \qquad\qquad\qquad\qquad\qquad\qquad\qquad\quad (1.1)\\
\ell(y)&=&-y''(x)+Q(x)y(x) \qquad\qquad\qquad\quad\qquad\quad\qquad\>(1.2)
\end{eqnarray*}

\noindent with the same boundary conditions

$$y'(0)=y'(\pi)=0 \qquad\qquad\qquad $$

\noindent  in $\>H_{1}\>$. We assume that the operator function $\>Q(x)\>$ in $\>(1.2)\>$
 verifies the following conditions:\\
\\
\noindent \textbf{(Q1)} \textit{Q(x)} has $\>(2k-2)$ th order weak derivatives $\>(k\geq2)\>$ and \\
$\>Q^{(i)}(x):H\rightarrow H \quad (i=0,1,2,\cdots,2k-2)\>$ are self-adjoint kernel operators for
 every $\>x\in[0,\pi]$.\\
\\
\noindent \textbf{(Q2)} $\> \parallel Q \parallel < \frac{1}{2}$.\\
\\
\noindent \textbf{(Q3)} There is an orthonormal basis $\>\{\varphi_{n}\}^{\infty}_{n=1}\>$ of $\>H\>$
such that

$$\sum\limits_{n=1}^{\infty}\parallel Q(x)\varphi_{n}\parallel <\infty.$$

\newpage
\noindent \textbf{(Q4)} The functions $\>\parallel Q^{(i)}(x)\parallel_{\sigma_{1}(H)} \quad (i=0,1,2,\cdots,2k-2)\>$
are both bounded and measurable in $\>[0,\pi]$.\\
\\
\noindent Here, $\>\sigma_{1}(H)\>$ denotes Banach space consisting of kernel operators from
$\>H\>$ to $\>H$.

\noindent Spectrum of the operator $\>L_{0}\>$ is the set $\>\{m^{2}\}_{m=0}^{\infty}$.
Each point of this set is an eigenvalue of $\>L_{0}\>$ which has
infinite multiplicity \cite{Hl}.

\noindent The eigenvectors corresponding to $\>m^{2}\>$ are the form:

$$\psi_{mn}=K_{m}\varphi_{n}\cos mx \qquad (m=0,1,2,\cdots \>;\>n=1,2,\cdots)\eqno(1.3)$$

\noindent where
\begin{displaymath}
K_m=\left\{\begin{array}{ll}
\frac{1}{\sqrt{\pi}}\>, &m=0\\
\\
\sqrt{\frac{2}{\pi}}\>, &m=1,2,\cdots
\end{array} \right.
\end{displaymath}

\noindent The research on the regularized trace of scalar differential
operators is studied for the first time by  $\cite{Ge}$.
In many other works such as $\cite{Di}$, $\cite{Ha}$ regularized traces of various
scalar differential operators have been investigated. The list of
the works on this subject is given in \cite{Le} and \cite{Fu}.
In the references $\cite{Hl}$, $\cite{Ma}$, $\cite{Ba}$ the regularized trace formulas
of
Sturm-Liouville operators with operator coefficient have been found. In \cite{Az}, the authors find
a trace formula for a self-adjoint operator in the following form:

\begin{eqnarray*}
\sum\limits_{m=1}^{\infty} \biggl\{\sum\limits_{n=1}^{\infty}\bigl( \lambda_{mn}^{2}&-& m^{4}\bigr)-\frac{2m^{2}}{\pi}\int\limits_{0}^{\pi}
\rm{tr}Q(x)\ud x-C\biggr\}\\
\\
&=&\frac{C}{2}+\frac{1}{8}\Bigl[\rm{tr}Q''(0)+\rm{tr}Q''(\pi)\Bigr]-\frac{1}{4}\Bigl[\rm{tr}Q^{2}(0)+\rm{tr}Q^{2}(\pi)\Bigr]
\end{eqnarray*}

\noindent where
$$ C=\frac{1}{2\pi}\Bigl[\rm{tr}Q'(0)-\rm{tr}Q'(\pi)\Bigr]+\frac{1}{2\pi}\int\limits_{0}^{\pi}\rm{tr}Q^{2}(x)\ud x+\frac{1}{2\pi^{2}}\>\rm{tr}
\biggl(\int\limits_{0}^{\pi}Q(x)\ud x\biggr)^{2}.$$

\noindent In the present paper the trace formula
\newcommand{\Res} {\mathrm {Res}}
\begin{eqnarray*}
\sum\limits_{m=0}^{\infty} \bigg\{\sum\limits_{n=1}^{\infty}\Big( \lambda_{mn}^{k}-m^{2k}\Big)&-&k\sum\limits_{j=2}^{2k+2}(-1)^{j}\>j^{-1}
{\Res}_{\lambda=m^{2}} \bigg[\lambda^{k-1}\rm{tr}\Big(QR_{\lambda}^{0}\Big)^{j}\bigg]\bigg\}\\
\\
&=&(-1)^{k-1}\>k\>2^{-2k}\Bigl[\rm{tr}Q^{2k-2}(0)+\rm{tr}Q^{2k-2}(\pi)\Bigr]
\end{eqnarray*}
\noindent is established.

\section{Some Formulas About Eigenvalues}

 Let $\>R_{\lambda}^{0}\>$ and $\>R_{\lambda}\>$ be resolvents of the operator $\>L_{0}\>$ and $\>L$,
respectively. Since the operator $\>Q(x)\>$ satisfies the condition $\> (Q3)\>$ then
$\> QR_{\lambda}^{0}:H_{1}\rightarrow H_{1}\>$  is a kernel operator for every
$\>\lambda\notin \{m^2\}_{0}^{\infty}\>$ we can show this  by using  the study $\cite{Co}$.

\noindent In $\cite{Hl}$, it is proved that the spectrum of the operator $\>L\>$ is a subset
of the union of pairs of disjoint intervals

$$ [m^2-\frac{1}{2}\>,\> m^2+\frac{1}{2}] \qquad (m=0,1,2,\cdots).$$

\noindent Moreover  it is proved that $\>m^2 \quad(m=0,1,\cdots)\>$ can be an eigenvalue of $\>L\>$
which has finite or infinite multiplicity and

$$ \lim_{n\rightarrow \infty}\lambda_{mn}=m^2$$

\noindent where $\>\{\lambda_{mn}\}_{n=1}^{\infty}\>$ are eigenvalues of $\>L\>$ belonging to the interval
$\>[m^2-\frac{1}{2}\>,\> m^2+\frac{1}{2}]\>$ in $\cite{Hl}$.

\noindent Each point belonging to $\> [m^2-\frac{1}{2}\>,\> m^2+\frac{1}{2}]$, different from
$\>m^2\>$ of the spectrum of $\>L\>$ is an isolated eigenvalue which has finite multiplicity. we can show this similar to $\cite{Yo}$

 Since $\>QR_{\lambda}^{0}\>$ is a kernel operator and
$$ R_{\lambda}-R_{\lambda}^{0}=-R_{\lambda}QR_{\lambda}^{0} \eqno(2.1)$$
\noindent we have $\>R_{\lambda}-R_{\lambda}^{0}\in \sigma_{1}(H_1)\>$ for each $\>\lambda\>$
which belongs to the resolvent set of $\>L\>$. In this case from $\cite{Co}$ we obtain the formula

$$ \rm{tr}\Big(R_{\lambda}-R_{\lambda}^{0}\Big)=\sum\limits_{m=0}^{\infty}\sum\limits_{n=1}^{\infty}
\Bigl(\frac{1}{\lambda_{mn}-\lambda}-\frac{1}{m^{2}-\lambda}\Bigr).$$

\noindent If we multiply the last equality by $\>\frac{\lambda^{k}}{2\pi i}\>$ and integrate on the circle \\
$\>|\lambda|=b_{p}=p^{2}+p \quad (p\geq1)\>$, then we have

$$ \frac{1}{2\pi i}\int\limits_{|\lambda|=b_{p}}\lambda^{k} \rm{tr}\big(R_{\lambda}-R_{\lambda}^{0}\big)\ud \lambda =\sum\limits_{m=0}^{p}
\sum\limits_{n=1}^{\infty} \Bigl(m^{2k}-\lambda_{mn}^{k}\Bigr). \eqno(2.2)$$

\noindent By using the formula (2.1) we obtain

$$ R_{\lambda}-R_{\lambda}^{0}=\sum\limits_{j=1}^{N}(-1)^{j}\>R_{\lambda}^{0} \bigl(QR_{\lambda}^{0}\bigr)^{j}+(-1)^{N+1}\>R_{\lambda}
\bigl(QR_{\lambda}^{0}\bigr)^{N+1}$$

\noindent where $\>N\geq 1\>$ is any natural number.

\noindent Using the last equality we can rewrite (2.2) as follows

$$\sum\limits_{m=0}^{p} \sum\limits_{n=1}^{\infty} \Bigl(\lambda_{mn}^{k}-m^{2k}\Bigr)=\sum\limits_{j=1}^{N} M_{pj}+M_{p}^{(N)}.\eqno(2.3)$$

\noindent Here,
$$ M_{pj}=\frac{(-1)^{j+1}}{2\pi i}\int\limits_{|\lambda|=b_{p}}\lambda^{k} \rm{tr}\Big[R_{\lambda}^{0}\Big(Q R_{\lambda}^{0}\Big)^{j}\Big]\ud
\lambda \eqno(2.4)$$

$$ M_{p}^{(N)}=\frac{(-1)^{N}}{2\pi i}\int\limits_{|\lambda|=b_{p}}\lambda^{k} \rm{tr}\Big[R_{\lambda}\Big(Q R_{\lambda}^{0}\Big)^{N+1}\Big]\ud
\lambda. \eqno(2.5)$$

\newtheorem{theorem}{Theorem}[section]
\newtheorem{lemma}[theorem]{Lemma}

\begin{theorem} If the operator function $\>Q(x)\>$ satisfies condition  $\> (Q3)\>$ then

$$ M_{pj}=\frac{(-1)^{j}k}{2\pi ij}\int\limits_{|\lambda|=b_{p}}\lambda^{k-1} \rm{tr}\Bigl(QR_{\lambda}^{0}\Bigr)^{j}\ud \lambda .\eqno(2.6)$$

\end{theorem}

\noindent {\bf Proof: } We can show that the operator function $\>QR_{\lambda}^{0}\>$ is analytic with
respect to the norm in the space $\>\sigma_{1}(H_{1})\>$ in domain $\>\rho(L_{0})= C\backslash \sigma(L_{0})\>$ and

$$ \rm{tr}\bigg\{\Big[\big(QR_{\lambda}^{0}\big)^{j}\Big]'\bigg\}=j\rm{tr}\Big[\big(QR_{\lambda}^{0}\big)'\big(QR_{\lambda}^{0}\big)^{j-1}\Big].
\eqno(2.7)$$

\noindent If we consider the equation $\>\big(QR_{\lambda}^{0}\big)'=\big(QR_{\lambda}^{0}\big)^{2}\>$ then
we can write the formula $\>(2.7)\>$ as follows

$$ \rm{tr}\bigg\{\Big[\big(QR_{\lambda}^{0}\big)^{j}\Big]'\bigg\}=j\rm{tr}\Big[R_{\lambda}^{0}\big(QR_{\lambda}^{0}\big)^{j}\Big]. \eqno(2.8)$$

\noindent If we write the formula $\>(2.8)\>$ in $\>(2.4)\>$ then we have

$$ M_{pj}=\frac{(-1)^{j+1}}{2\pi ij}\int\limits_{|\lambda|=b_{p}}\lambda^{k} \rm{tr}\bigg\{\Big[\big(Q R_{\lambda}^{0})^{j}\Big]'\bigg\}
\ud \lambda. $$

\noindent Here we find

\begin{eqnarray*}
M_{pj}&=&\frac{(-1)^{j+1}}{2\pi ij}\int\limits_{|\lambda|=b_{p}} \rm{tr}\bigg\{\Big[\lambda^{k} \big(Q R_{\lambda}^{0}\big)^{j}\Big]'
-k\lambda^{k-1}\big(Q R_{\lambda}^{0}\big)^{j}\bigg\}\ud \lambda\\
\\
&=&\frac{(-1)^{j}}{2\pi ij}\int\limits_{|\lambda|=b_{p}} k\lambda^{k-1} \rm{tr}\big(Q R_{\lambda}^{0}\big)^{j}\ud \lambda\\
 \\
&+&\frac{(-1)^{j+1}}{2\pi ij}\int\limits_{|\lambda|=b_{p}} \rm{tr}\bigg\{\Big[\lambda^{k} \big(Q R_{\lambda}^{0}\big)^{j}\Big]'\bigg\}\ud \lambda
. \qquad \qquad\qquad\qquad\qquad\quad\quad\quad(2.9)
\end{eqnarray*}

\noindent We can show that

$$ \rm{tr}\bigg\{\Big[\lambda^{k}\big(Q R_{\lambda}^{0})^{j}\Big]'\bigg\}= \bigg\{\rm{tr}\Big[\lambda^{k} \big(Q
R_{\lambda}^{0}\big)^{j}\Big]\bigg\}'.$$

\noindent Therefore we have

$$\int\limits_{|\lambda|=b_{p}} \rm{tr}\bigg\{\Big[\lambda^{k} \big(Q R_{\lambda}^{0}\big)^{j}\Big]'\bigg\}\ud
\lambda=\int\limits_{|\lambda|=b_{p}} \bigg\{\rm{tr}\Big[\lambda^{k} \big(Q R_{\lambda}^{0}\big)^{j}\Big]\bigg\}'\ud \lambda.
\qquad\eqno(2.10)$$

\noindent We can write the integral on the write hand side of the equality $\>(2.10)\>$  as follows

\begin{eqnarray*}
\int\limits_{|\lambda|=b_{p}} \bigg\{\rm{tr}\Big[\lambda^{k} \big(Q R_{\lambda}^{0}\big)^{j}\Big]\bigg\}'\ud \lambda &=&
\int\limits_{|\lambda|=b_{p}\atop{Im\lambda\geq 0}} \bigg\{\rm{tr}\Big[\lambda^{k} \big(Q R_{\lambda}^{0}\big)^{j}\Big]\bigg\}'\ud \lambda\
\\
&+& \int\limits_{|\lambda|=b_{p}\atop{Im\lambda\leq 0}} \bigg\{\rm{tr}\Big[\lambda^{k} \big(Q R_{\lambda}^{0}\big)^{j}\Big]\bigg\}'\ud \lambda.
\qquad\qquad\quad (2.11)
\end{eqnarray*}

\noindent Let $\>\varepsilon_{0}\>$ be a constant which satisfies the condition $\> 0<\varepsilon_{0}< b_{p}-(p^{2}+p)$.

\noindent Since the function $\>tr \biggl[\lambda^{k} \big(Q R_{\lambda}^{0}\big)^{j}\biggr]\>$ is analytic in simple connected domains
\begin{eqnarray*}
G_{1}&=&\big\{\lambda\in C: b_{p}-\varepsilon_{0}<|\lambda|<b_{p}+\varepsilon_{0},\> Im\lambda > -\varepsilon_{0}\big\}\\
G_{2}&=&\big\{\lambda\in C: b_{p}-\varepsilon_{0}<|\lambda|<b_{p}+\varepsilon_{0},\> Im\lambda <-\varepsilon_{0} \bigr\}
\end{eqnarray*}

\noindent and
\begin{eqnarray*}
\bigl\{\lambda\in C:\> |\lambda|=b_{p},\> Im\lambda \geq 0\bigr\}\subset G_{1}\\
\bigl\{\lambda\in C:\> |\lambda|=b_{p},\> Im\lambda \leq 0\bigr\}\subset G_{2}
\end{eqnarray*}

\noindent then we rewrite equality $\>(2.11)\>$ 

\begin{eqnarray*}
\int\limits_{|\lambda|=b_{p}} \bigg\{\rm{tr}\Big[\lambda^{k} \Big(Q R_{\lambda}^{0}\Big)^{j}\Big]\bigg\}'\>\ud \lambda &=& \rm{tr}\Big[-b_{p}
\Big(QR_{-b_{p}}^{0} \Big)^{j}\Big]-\rm{tr}\Big[b_{p}\big(QR_{b_{p}}^{0}\big)^{j}\Big]\
\\
&+& \rm{tr}\Big[b_{p}\big(QR_{b_{p}}^{0} \big)^{j}\Big]-\rm{tr}\Big[-b_{p}\big(QR_{-b_{p}}^{0} \big)^{j}\Big]\\
\\
&=&0. \qquad\qquad\qquad\qquad\qquad\qquad\qquad(2.12)
\end{eqnarray*}

\noindent From (2.9), (2.10) and (2.12) we find

$$ M_{pj}=\frac{(-1)^{j}k}{2\pi ij}\int\limits_{|\lambda|=b_{p}} \lambda^{k-1}\rm{tr}\bigl(Q R_{\lambda}^{0}\bigr)^{j}\>\ud \lambda .$$

\section{Regularized Trace Formula of the Operator L}

In this section, we will obtain a regularized trace formula for the operator $\>L$.

 According to Theorem 2.1,

$$ M_{p1}=\frac{-k}{2\pi i}\int\limits_{|\lambda|=b_{p}} \lambda^{k-1}\rm{tr}\bigl(Q R_{\lambda}^{0}\bigr)\ud \lambda . \eqno(3.1)$$

\noindent Since $\>\{\psi_{mn}^{0} \}_{m=0,\>n=1}^{\infty\quad\infty}\>$ is an orthonormal basis of the space $\>H_{1}\>$, then we obtain

\begin{eqnarray*}
M_{p1}&=&\frac{-k}{2\pi i}\int\limits_{|\lambda|=b_{p}} {\lambda}^{k-1}\sum\limits_{m=0}^{\infty}\sum\limits_{n=1}^{\infty}\Bigl(Q
{R_{\lambda}}^{0}\psi_{mn},\psi_{mn}\Bigr)\ud \lambda\\
\\
&=& k\sum\limits_{m=0}^{\infty}\sum\limits_{n=1}^{\infty}\Bigl(Q \psi_{mn},\psi_{mn}\Bigr)\>\frac{1}{2\pi\>i} \int\limits_{|\lambda|=b_{p}}
\frac{{\lambda}^{k-1}}{\lambda-m^{2}}\>\ud \lambda \\
\\
 &=& k\sum\limits_{m=0}^{p}\sum\limits_{n=1}^{\infty} m^{2k-2}\Bigl(Q \psi_{mn},\psi_{mn}\Bigr)\\
\\
&=& k\sum\limits_{m=1}^{p}\sum\limits_{n=1}^{\infty} m^{2k-2} \, \frac{2}{\pi}\int\limits_{0}^{\pi}\Bigl(Q(x)\varphi_{n},\varphi_{n}\Bigr)_{H}\>
\cos^{2}mx\>\ud x\\
\\
&=& \frac{k}{\pi}\sum\limits_{m=1}^{p}\sum\limits_{n=1}^{\infty} m^{2k-2}\int\limits_{0}^{\pi}\Bigl(Q(x)\varphi_{n},\varphi_{n}\Bigr)_{H}\>
\Bigl(1+\cos2mx\Bigr)\ud x\\
\\
&=&\frac{k}{\pi}\sum\limits_{m=1}^{p} m^{2k-2}\int\limits_{0}^{\pi} \rm{tr}Q(x)\ud x +\frac{k}{\pi}\sum\limits_{m=1}^{p}
m^{2k-2}\int\limits_{0}^{\pi} \rm{tr}Q(x)\cos2mx\>\ud x \>\qquad(3.2)
\end{eqnarray*}

\noindent for $\>k\geq2$.
\noindent Let us evaluate the following integral

{\setlength\arraycolsep{2pt}
\begin{eqnarray*}
\int\limits_{0}^{\pi}&\rm{tr}&Q(x)\cos2mx\ud x=\frac{1}{2m}\sin 2mx\>\rm{tr}Q(x)\Big|_{0}^{\pi} -\frac{1}{2m}\int\limits_{0}^{\pi}
\rm{tr}Q'(x)\sin 2mx\>\ud x {}\\
\\
&&{}=-\frac{1}{2m}\int\limits_{0}^{\pi} \rm{tr}Q'(x)\sin 2mx\>\ud x{}\\
\\
&&{}=\frac{1}{2m}\biggl[\frac{1}{2m}\cos 2mx\>\rm{tr}Q'(x)\Big|_{0}^{\pi}-\frac{1}{2m}\int\limits_{0}^{\pi} \rm{tr}Q''(x)\cos 2mx\>\ud
x\biggr]{}\\
\\
&&{}=\frac{1}{4m^2}\biggl(\rm{tr}Q'(\pi)-\rm{tr}Q'(0)\biggr)-\frac{1}{4m^2}\int\limits_{0}^{\pi} \rm{tr}Q''(x)\cos 2mx\>\ud x{}\\
\\
&&{}=\frac{1}{4m^2}\biggl(\rm{tr}Q'(\pi)-\rm{tr}Q'(0)\biggr)-\frac{1}{4m^2}\biggl[\frac{1}{2m}\sin 2mx\>\rm{tr}Q''(x)\Big|_{0}^{\pi}{}\\
\\
&&{}-\frac{1}{2m}\int\limits_{0}^{\pi} \rm{tr}Q'''(x)\sin 2mx\>\ud x\biggr]{}\\
\\
&&{}=\frac{1}{4m^2}\biggl(\rm{tr}Q'(\pi)-\rm{tr}Q'(0)\biggr)+\frac{1}{8m^3}\int\limits_{0}^{\pi} \rm{tr}Q'''(x)\sin 2mx\>\ud x{}\\
\\
&&{}=\frac{1}{4m^2}\biggl(\rm{tr}Q'(\pi)-\rm{tr}Q'(0)\biggr)-\frac{1}{8m^3}\bigg[\frac{1}{2m}\cos 2mx\>\rm{tr}Q'''(x)\Big|_{0}^{\pi}{}\\
\\
&&{}-\frac{1}{2m}\int\limits_{0}^{\pi} \rm{tr}Q^{\i v}(x)\cos 2mx\>\ud x\>\bigg]{} \\
\\
&&{}=\frac{1}{4m^2}\biggl(\rm{tr}Q'(\pi)-\rm{tr}Q'(0)\biggr)-\frac{1}{16m^4}\Big(\rm{tr}Q'''(\pi)-\rm{tr}Q'''(0)\Big){}\\
\\
&&{}+\frac{1}{16m^4}\int\limits_{0}^{\pi}\rm{tr}Q^{\i v}(x)\cos 2mx\>\ud x{}\\
\\
&&{}=\frac{1}{4m^2}\Big(\rm{tr}Q'(\pi)-\rm{tr}Q'(0)\Big)-\frac{1}{16m^4}\Big(\rm{tr}Q'''(\pi)-\rm{tr}Q'''(0)\Big){}\\
\\
&&{}+\frac{1}{16m^4}\bigg[\frac{1}{2m}\sin 2mx\>\rm{tr}Q^{\i v}(x)\Big|_{0}^{\pi}-\frac{1}{2m}\int\limits_{0}^{\pi}\rm{tr}Q^{v}(x)\sin 2mx\>\ud
x\bigg]{}\\
\\
&&{}=\frac{1}{4m^2}\Bigl(\rm{tr}Q'(\pi)-\rm{tr}Q'(0)\Bigr)-\frac{1}{16m^4}\Big(\rm{tr}Q'''(\pi)-\rm{tr}Q'''(0)\Big)\\
\\
&&{}+\frac{1}{32m^5}\Bigl[\frac{1}{2m}\cos 2mx\>\rm{tr}Q^{v}(x)\Big|_{0}^{\pi}-\frac{1}{2m}\int\limits_{0}^{\pi}\rm{tr}Q^{v\i}(x)\cos 2mx\>\ud
x\Bigr]{}\\
\\
&&{}=\frac{1}{4m^2}\Big(\rm{tr}Q'(\pi)-\rm{tr}Q'(0)\Big)-\frac{1}{16m^4}\Big(\rm{tr}Q'''(\pi)-\rm{tr}Q'''(0)\Big){}\\
\\
&&{}+\frac{1}{64m^6}\Big(\rm{tr}Q^{v}(\pi)-\rm{tr}Q^{v}(0)\Big){}
\\
&&{}-\frac{1}{64m^6}\int\limits_{0}^{\pi}\rm{tr}Q^{v\i}(x)\cos2mx\>\ud x{}\\
\\
&&{} = \cdots \\
\\
&&{}=\sum\limits_{i=2}^{k}\frac{(-1)^{i}}{(2m)^{2i-2}}\Bigl(\rm{tr}Q^{(2i-3)}(\pi)-\rm{tr}Q^{(2i-3)}(0)\Bigr){}\\
\\
&&{}+\frac{(-1)^{k-1}}{(2m)^{2k-2}} \int\limits_{0}^{\pi} \rm{tr}Q^{(2k-2)}(x)\cos2mx\>\ud x \qquad\qquad\qquad\qquad\qquad (3.3)
\end{eqnarray*}}

\noindent From (3.2) and (3.3) we have

\begin{eqnarray*}
M_{p1}&=&\frac{k}{\pi}\sum\limits_{m=1}^{p}m^{2k-2}\int\limits_{0}^{\pi} \rm{tr}Q(x)\>\ud x\\
\\
&+& \frac{k}{\pi}\sum\limits_{m=1}^{p}\sum\limits_{i=2}^{k} (-1)^{i}2^{2-2i}m^{2k-2i}\Bigl(\rm{tr}Q^{(2i-3)}(\pi)-\rm{tr}Q^{(2i-3)}(0)\Bigr)\\
\\
&-&\frac{k}{\pi}\sum\limits_{m=1}^{p}(-1)^{k}\>2^{2-2k}\int\limits_{0}^{\pi} \rm{tr}Q^{(2k-2)}(x)\cos2mx\>\ud x \qquad\qquad\qquad\qquad(3.4)
\end{eqnarray*}

\noindent From (2.3) and (3.4) we have

\begin{eqnarray*}
\sum\limits_{m=0}^{p} \bigg[\sum\limits_{n=1}^{\infty} \Bigl(\lambda_{mn}^{k}&-&m^{2k}\Bigr)-\>\sum\limits_{j=2}^{N} M_{pj}-k{\pi}^{-1}m^{2k-2}
\int\limits_{0}^{\pi} \rm{tr}Q(x)\>\ud x\\
&-&k{\pi}^{-1}\sum\limits_{i=2}^{k}(-1)^{i}\>2^{2-2i}m^{2k-2i}\Big(\rm{tr}Q^{(2i-3)}(\pi)-\rm{tr}Q^{(2i-3)}(0)\Big)\bigg]\\
&=&-k{\pi}^{-1}\sum\limits_{m=1}^{p}(-1)^{k}2^{2-2k}\int\limits_{0}^{\pi} \rm{tr}Q^{(2k-2)}(x)\cos2mx\ud x\\
&+& M_{p}^{(N)} \qquad\qquad\qquad\qquad\qquad\qquad\qquad\qquad\qquad\qquad\quad(3.5)
\end{eqnarray*}

\noindent Similar to $\cite{Yo}$, it can be shown that

$$ \mid M_{p}^{N}\mid \leq \textrm{cons}t\>p^{1+2k-N}.\eqno(3.6)$$

\noindent From $\>(3.6)\>$ we obtain

$$\lim_{p\rightarrow\infty}{M_{p}}^{(2k+2)}=0. \eqno(3.7)$$

\noindent On the other hand from  $\cite{Ad}$ we have

{\setlength\arraycolsep{2pt}
\begin{eqnarray*}
\frac{1}{\pi}\sum\limits_{m=1}^{\infty}\int\limits_{0}^{\pi}\rm{tr}Q^{(2k-2)}(x)\cos2mx\>\ud x&=&\frac{1}{4}\Bigl[\rm{tr}Q^{(2k-2)}(0)\>+\>
\rm{tr}Q^{(2k-2)}(\pi)\Bigr]\\
\\
&-&\frac{1}{2\pi}\Bigl[\rm{tr}Q^{(2k-3)}(\pi) -\rm{tr}Q^{(2k-3)}(0)\Bigr].  \quad\>(3.8)
\end{eqnarray*}}

\noindent Moreover the expression $\>(2.6)\>$ can be written as

{\setlength\arraycolsep{2pt}
\begin{eqnarray*}
 M_{pj}&=&(-1)^{j}j^{-1}k\>\frac{1}{2\pi i}\int\limits_{|\lambda|=b_{p}}\lambda^{k-1} \rm{tr}\Bigl(QR_{\lambda}^{0}\Bigr)^{j}\ud \lambda {}\\
 \\
 &=&(-1)^{j}j^{-1}k\sum\limits_{m=0}^{p}\rm{Res}_{\lambda=m^{2}}\biggl[\lambda^{k-1}tr\big(QR_{\lambda}^{0}\big)^{j}\biggr]
 \qquad\qquad\qquad\qquad\>\>(3.9)
 \end{eqnarray*}}

\noindent If we consider the formulas $\>(3.7),\>(3.8),\>(3.9)\>$ and we take limit in $\>(3.5)\>$ as $\>p \rightarrow \infty \>$ we obtain
$ k\emph{th}$ order regularized trace formula of the operator $\>L$.

{\setlength\arraycolsep{2pt}
\begin{eqnarray*}
\sum\limits_{m=0}^{\infty}\biggl\{\sum\limits_{n=1}^{\infty}\Big({\lambda}_{mn}^{k}&-&{m}^{2k}\Big)-k\sum\limits_{j=2}^{2k+2}(-1)^{j}j^{-1}
Res_{\lambda=m^{2}}\bigg[{\lambda}^{k-1}\mbox{tr}\Big(QR_{\lambda}^{0}\Big)^{j}\bigg]\\
\\
&-&k{\pi}^{-1}m^{2k-2}\int\limits_{0}^{\pi}\mbox{tr}Q(x)\>\ud x\\
\\
&-& k{\pi}^{-1}\sum\limits_{i=2}^{k} {(-1)}^{i}\>2^{2-2i}\>m^{2k-2i}\bigg(\mbox{tr}Q^{(2i-3)}(\pi)- \mbox{tr}Q^{(2i-3)}(0)\bigg)\biggr\}\\
\\
&=&{(-1)}^{k-1}\>2^{-2k}\Big[\mbox{tr}Q^{(2k-2)}(0)+ \mbox{tr}Q^{(2k-2)}(\pi)\Big]\\
\\
&+&{(-1)}^{k}{k}\>{\pi}^{-1}2^{1-2k}\Big[\mbox{tr}Q^{(2k-3)}(\pi)- \mbox{tr}Q^{(2k-3)}(0)\Big]
\end{eqnarray*}}

\noindent as $\>p\longrightarrow \infty\>$.

\noindent  We can write the k\emph{th} regularized trace formula of $\>L\>$ as follows

{\setlength\arraycolsep{2pt}
\begin{eqnarray*}
\sum\limits_{m=0}^{\infty}\biggl\{\sum\limits_{n=1}^{\infty}\Big({\lambda}_{mn}^{k}&-&{m}^{2k}\Big)- k\sum\limits_{j=2}^{2k+2}(-1)^{j}j^{-1}
Res_{\lambda=m^{2}}\biggl[\lambda^{k-1}\mbox{tr}\big(QR_{\lambda}^{0}\big)^{j}\bigg]{}\\
\\
&-&\>k{\pi}^{-1}m^{2k-2}\int\limits_{0}^{\pi}\mbox{tr}Q(x)\>dx\>-\>4k{\pi}^{-1}m^{2k}\sum\limits_{i=2}^{k} m^{-2i}a_{i}\biggr\}\\
\\
&=&{(-1)}^{k-1}k2^{-2k}\Big[\mbox{tr}Q^{(2k-2)}(0)\>+\> \mbox{tr}Q^{(2k-2)}(\pi)\Big]{}\\
\\
&+& 2{\pi}^{-1}ka_{k} \qquad\qquad\qquad\qquad\qquad\qquad\qquad\quad(3.10)
\end{eqnarray*}}

\noindent where $\>a_{i}={(-1)}^{i}2^{-2i}\Big[\mbox{tr} Q^{(2i-3)}(\pi)- \mbox{tr}Q^{(2i-3)}(0)\Big] \quad (i=2,3,\cdots,k)$.

\noindent In addition if $\>Q(x)\>$ satisfies the following conditions

\noindent $\textbf{(Q5)}\> Q^{(2i-3)}(0)=Q^{(2i-3)}(\pi)=0 \quad (2\leq i\leq k)$

\noindent$ \textbf{(Q6)}\> \int\limits_{0}^{\pi}\mbox{tr}Q(x)\>\ud x=0$

\noindent then we can write formula $\>(3.10)\>$

{\setlength\arraycolsep{2pt}
\begin{eqnarray*}
\sum\limits_{m=0}^{\infty}\biggl\{\sum\limits_{n=1}^{\infty}\Big({\lambda}_{mn}^{k}-{m}^{2k}\Big)&-& k\sum\limits_{j=2}^{2k+2}(-1)^{j}j^{-1}
Res_{\lambda=m^{2}}\Bigl[\lambda^{k-1}\mbox{tr}\big(QR_{\lambda}^{0}\big)^{j}\Big]\biggr\}{}\\
\\
&=&(-1)^{k-1}k\>2^{-2k}\Big[\mbox{tr}Q^{(2k-2)}(0)\>+\> \mbox{tr}Q^{(2k-2)}(\pi)\Big]
\end{eqnarray*}}

\end{document}